\documentclass[aps,pre,twocolumn]{revtex4-1}

\usepackage[dvips]{graphicx}
\usepackage{amssymb,amsfonts,amsmath}
\usepackage{color}
\usepackage{colortbl}
\usepackage{wasysym}
\usepackage{ifsym}

\begin{document}
\title{Signature stability analysis for networks of coupled dynamical systems with Hermitian Jacobian}
\author{Anne--Ly Do,$^1$ Stefano Boccaletti,$^2$ Jeremias Epperlein,$^3$ Stefan Siegmund,$^3$ and Thilo Gross $^4$}
\email[]{ly@mpipks-dresden.mpg.de}
\affiliation{\mbox{$^1$ Max-Planck-Institute for the Physics of Complex Systems, Dresden, Germany}\\
\mbox{$^2$ Technical University of Madrid, Center for Biomedical Technologies, Madrid, Spain}\\
\mbox{$^3$ Dresden University of Technology, Center for Dynamics, Dresden, Germany}\\
\mbox{$^4$ University of Bristol, Merchant Venturers School of Engineering, Bristol, United Kingdom}}
\date{\today}

\begin{abstract}
The central theme of complex systems research is understanding the emergent macroscopic properties of a system from the interplay of its microscopic constituents.
Here, we ask what conditions a complex network of microscopic dynamical units has to meet to permit stationary macroscopic dynamics, such as stable equilibria or phase-locked states. 
We present an analytical approach which is based on a graphical notation that allows rewriting Jacobi's signature criterion in an interpretable form.
The derived conditions pertain to topological structures on all scales, ranging from individual nodes to the interaction network as a whole.
Our approach can be applied to many systems of symmetrically coupled units. 
For the purpose of illustration, we consider the example of synchronization, specifically the (heterogeneous) Kuramoto model and an adaptive variant.
Moreover, we discuss how the graphical notation can be employed to study isospectrality in Hermitian matrices.  
The results complete and extend the previous analysis of Do et al. [Phys. Rev. Lett. \textbf{108}, 194102 (2012)]. 
\end{abstract}

\maketitle

Characterizing the behavior of complex systems and discovering the critical boundaries in parameter space at which qualitative changes occur are of central interest in statistical physics.
Many studies have demonstrated that, in specific models, this can be achieved by the instruments of nonlinear dynamics \cite{Laradji,Atay,Valladares,Eckhardt}. 
However, applying these tools to large heterogeneous systems, such as complex networks of dynamical units poses significant challenges \cite{Barabasi,Newman,Newman2006,Dorogovtsev}. 
In particular, there is no simple way to capture the structural aspects of such systems, which are known to crucially influence the emergent dynamics \cite{Boccaletti,Fortunato}.
  
The interplay of dynamics and structure in networks has extensively been studied using the example of synchronization \cite{Kurths,Boccaletti,Arenas}. 
Systems of coupled phase oscillators, originally proposed by Kuramoto \cite{Kuramoto,Acebron}, have become a paradigmatic example. 
The exploration of such systems revealed that the synchronizability of a network depends on global topological measures such as the clustering coefficient, the diameter, and the degree or weight distribution \cite{Mirollo,Wu,Pecora,Kawamura,Chavez,Nishikawa2006,Lodato}, but also on local and mesoscale details \cite{Arenas,Nishikawa,Mori}.  
This shows that a conceptual understanding of collective phenomena on dynamical networks can only be achieved if structural properties on all scales are taken into account. 

In a recent paper, we proposed a nonlinear dynamics approach, which allows to study the dependence of dynamics on structure in networks of symmetrically coupled units \cite{Do,Epperlein}. 
We showed that Jacobi's signature criterion (JSC) can be used to determine necessary conditions an interaction network has to obey in order to support a given type of macroscopic behavior. 
These conditions pertain to subgraphs on all scales, from single nodes to the entire network. 
In this paper, we complete and extend the analysis of \cite{Do}. 
We present the full derivation of the proposed approach and discuss a number of new applications including an adaptive Kuramoto model, in which the topology of the coupling network coevolves with the local dynamics of the phase oscillators \cite{ZhouKurths,Almendral}.
  
\section{Stability in networks of phase-oscillators}  
The approach outlined below is applicable to all dynamical systems with Hermitian Jacobian matrices, and in particular to all 
networks of symmetrically-coupled one-dimensional dynamical units. However, for the sake of illustration, we first consider   
specifically the example of $N$ coupled phase oscillators
\begin{equation}\label{ODE}
\dot{x_i}=\omega_i+\sum_{j\neq i} A_{ij}\sin(x_j-x_i)\ , \quad \forall i\in 1\ldots N	\ .
\end{equation}
Here, $x_i$ and $\omega_i$ denote the phase and the intrinsic frequency of node $i$, and ${\rm \bf A}\in \mathbb{R}^{N\times N}$ is the weight matrix of an undirected, weighted interaction network. 
Two oscillators $i,j$ are thus connected if $A_{ij}=A_{ji}\neq 0$.
Equation~\eqref{ODE} represents the so-called Kuramoto model, that is today considered to be a paradigm for the study of synchronization phenomena in coupled discrete systems \cite{Acebron}. 
The model is therefore used as the natural benchmark for comparative evaluations of performances of methods and tools. 

It is known that the Kuramoto-model can exhibit phase-locked states, which correspond to steady states of the governing system of equations. 
The local stability of such states is determined by the spectrum of the Jacobian matrix ${\rm \bf J} \in\mathbb{R}^{N\times N}$ defined by $J_{ik}=\partial\dot{x_{i}}/\partial x_{k}$. As ${\rm \bf J}$ is symmetric, its spectrum is real. The state under consideration is asymptotically stable if all eigenvalues are $\leq 0$.

Below, we analyze the spectrum of ${\rm \bf J}$ by means of Jacobi's signature criterion (JSC).
The JSC (also known as Sylvester criterion) states that a hermitian or symmetric matrix ${\rm \bf J}$ with rank $r=N-1$ has $r$ negative eigenvalues if and only if all principal minors of order $q\leq r$ have the sign of $(-1)^q$ \cite{Epperlein}.
Here, the principal minor of order $q$ is defined as $D_{q}:=\det\left(J_{ik}\right)$, $i,k=s_1,\ldots,s_q$.

Stability analysis by means of JSC is well-known in control theory \cite{LiaoYu} and has been applied to problems of different fields. 
For instance, the JSC has been used in \cite{Beckers} to analyze advection-diffusion equations, and in \cite{Soldatova} to determine stability boundaries in thermodynamic systems. 
Moreover, it has been applied to study the dynamics of offshore platforms \cite{CaiWuChen} and social systems \cite{DoRudolfGross}. 

Two obstacles that are encountered in the application of the JSC to large networked systems are the growth of
(a) the number of terms in each determinant, and (b) the number of determinants that have to be checked.  
Difficulty (a) can be overcome by means of a graphical notation, which allows to derive topological stability criteria. This notation is introduced in Sec.~\ref{Notation_sec}.  
Difficulty (b) can be evaded by considering necessary rather than sufficient conditions for stability: 
While the sufficient condition for stability necessitates the evaluation of all principal minors, a necessary condition for stability is already obtained by demanding \mbox{${\rm sgn}\left(D_q\right)=(-1)^q$} for some $q$. 

The necessary stability condition that is found by considering a principal minor of given order $q$ depends on the ordering of variables, i.e., the ordering of rows and columns in the Jacobian. 
By considering different orderings the number of conditions obtained for a given $q$ can therefore be increased \cite{DoRudolfGross}. 
To distinguish minors that are based on different orderings of the variables, we define $S=\left\{s_1,\ldots,s_q\right\}$ as a set of $\left| S\right|=q$ indices and $D_{\left| S\right|,S}$ as the determinant of the submatrix of $\rm \bf J$, which is spanned by the variables $x_{s_1},\ldots, x_{s_q}$.
Therewith, the conditions for stability read
\begin{equation}\label{stability_cond}
	{\rm sgn}\left(D_{\left| S\right|,S}\right)=(-1)^{\left| S\right|}, \quad \forall S,  \ \left| S\right|=1,\ldots, r.
\end{equation}

\section{Graphical notation}\label{Notation_sec}
Below, we address difficulty (a) mentioned above, i.e., the combinatorial explosion of terms that are needed to write out the conditions for increasing $\left| S\right|$. 
In the common notation, more than $700$ terms are necessary for expressing the minors of order $6$. 
For circumventing this problem we employ a graphical notation that captures basic intuition and allows for expressing the minors in a concise way \cite{Do}. 

The graphical notation relies on a topological reading of the minors. 
We interpret the Jacobian ${\rm \bf J}$ as the weight matrix of an undirected, weighted graph $\mathcal{G}$.
A Jacobian element $J_{ij}$ then corresponds to the weight of a link connecting nodes $i$ and $j$.
Moreover, we relate products of the Jacobian elements to subgraphs of $\mathcal{G}$ spanned by the respective links. 
For instance $J_{ij}J_{jk}$ is interpreted as path $i-j-k$. Similarly, $J_{ij}J_{jk}J_{ki}$ is interpreted as a closed path from $i$ to $j$ to $k$ and back to $i$. 

We can now express the minors of ${\rm \bf J}$ as sums over subgraphs of $\mathcal{G}$.
The Leibniz formula for determinants \cite{LA} implies that
(i) a minor $D_{|S|,S}$ is a sum over $|S|!$ elementary products $J_{i_1j_1}\cdot\ldots\cdot J_{i_{|S|}j_{|S|}}$; and that (ii) in each of these products every index $s_i \in S$ occurs exactly twice.
In the topological reading this translates to the following statements:
Because of property (i), each term of a minor $D_{|S|,S}$ corresponds to a subgraph with $|S|$ links.
Because of property (ii), these subgraphs can be decomposed in cycles of $\mathcal{G}$:
Every index $s_i\in S$ occurs either with multiplicity two on a diagonal element of ${\rm \bf J}$, or, with multiplicity one, on two off-diagonal elements of ${\rm \bf J}$.
In the former case, the respective factor corresponds to a self-loop of $\mathcal{G}$, i.e., to a cycle of length $n=1$; in the latter case, there is a set of factors $J_{ij}$, $i\neq j$, corresponding to a closed path of links, i.e., a cycle of length $n>1$.

In summary, the index structure of every term of a minor $D_{|S|,S}$ can be expressed by a combination of symbols denoting cycles of $\mathcal{G}$ of a given length $n$.
The idea is now to supplement the basis of symbols with a summation convention;
This convention is designed such that all algebraic terms that are structurally identical and only differ by index permutations can be captured in one symbolic term, which drastically reduces the complexity of the expressions. 
 
\begin{figure}
\includegraphics[width=0.5\textwidth]{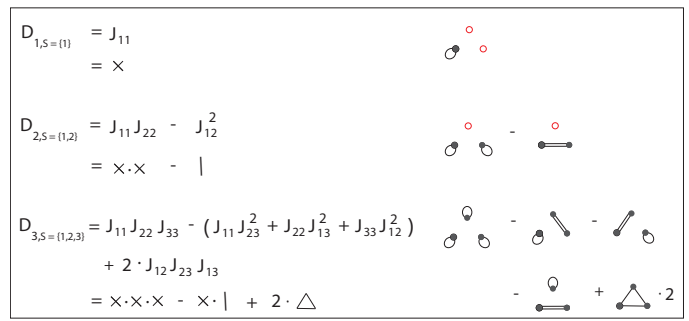} 
\vspace{-0.4cm}
\caption{Example for the graphical notation. Shown are the minors of the matrix~\eqref{Ex} in algebraic, and graphical notation, and, for each term, the corresponding subgraph of a three-node graph $\mathcal{G}$. Here, as well as in the next figures, filled symbols correspond to nodes $\in S$, open symbols to nodes $\notin S$.}
\label{Examples}
\end{figure} 
 
Below, we use the following definitions: 
The basis of symbols is given by $\times, |, \triangle, \Square, \pentagon, \ldots$ denoting cycles of length $n=1,2,3,4,5,\ldots$.
The summation convention stipulates that in a minor $D_{\left| S\right|,S}$, every product of symbols denotes the sum over all non-equivalent possibilities to build the depicted subgraph with the nodes $\in S$ (cf. Fig.~1).
With these conventions the first $4$ principal minors can be written as
\begin{subequations}\label{det_n}
\begin{align}
D_{1,S}&=\times\\
D_{2,S}&=\times \cdot \times -| \\
D_{3,S}&=\times \cdot \times \cdot \times - \times \cdot  | +2\triangle\\
D_{4,S}&=\times \cdot \times\cdot \times\cdot \times -\times\cdot\times \cdot |+|\cdot|+2\times \cdot \triangle
- 2\Square\label{EqFig2}
\end{align}
\end{subequations}
More generally, 
\begin{equation}\label{formationrule}
	D_{\left| S\right|,S}=\sum\text{all combinations of symbols with}\;\sum n=\left| S\right|,
\end{equation}
where symbols with $n>2$ appear with a factor of $2$ that reflects the two possible orientations in which the corresponding subgraphs can be paced out.
Symbols with an even (odd) number of links carry a negative (positive) sign related to the sign of the respective index permutation in the Leibniz formula for determinants \cite{LA}.

An example for the graphical notation\index{graphical notation} is presented in Fig.~\ref{Examples}. 
The figure displays the three principal minors of the symmetric $3\times 3$ matrix 
\begin{equation}\label{Ex} 
{\rm \bf J}\!=\!\left(\begin{array}{ccc}
J_{11}&J_{12}&J_{13} \\
 J_{12}&J_{22}&J_{23} \\
 J_{13}&J_{23}&J_{33}\\ 
\end{array}\right)\ 
\end{equation}
in algebraic, and graphical notation. Moreover, it displays for each term the corresponding subgraph of a three-node graph $\mathcal{G}$.

\section{Zero row sum}
In many systems, including the standard Kuramoto model, fundamental conservation laws impose a zero-row-sum condition, such that $J_{ii}=-\sum_{j\neq i}J_{ij}$.
Using this relation, we can remove all occurrences of elements $J_{ii}$ from the Jacobian and its minors \cite{Do}. 
In the topological reading, this substitution changes the graph $\mathcal{G}$ by replacing a self-loop at a node $i$ by the negative sum over all links that connect to $i$.

The simplification of the minors due to the zero-row-sum condition can be understood using the example of the Eqs.~\eqref{det_n}. 
Replacing the self-loops, the first term of every minor $D_{\left| S\right|,S}$, $\times^{\left| S\right|}$, is $(-1)^{\left| S\right|}$ times the sum over all subgraphs meeting the following conditions: (i) the subgraph contains exactly $\left| S\right|$ links, and (ii) it can be drawn by starting each (undirected) link at a different node in $S$. 
By means of elementary combinatorics it can be verified that all other terms of $D_{\left| S\right|,S}$ cancel exactly those subgraphs in $\times^{\left| S\right|}$ that contain cycles. This enables us to express the minors in another way: 
Defining 
\begin{equation}\label{def}\begin{small}
\Phi_{\left| S\right|,S}=\sum\text{all acyclic subgraphs of $\mathcal{G}$ meeting (i) and (ii)}
\end{small}
\end{equation}
we can write
\begin{equation}\label{det_n_2}
	D_{\left| S\right|,S}=(-1)^{\left| S\right|}\Phi_{\left| S\right|,S} \ .
\end{equation}

We remark that for $\left| S\right|=N-1$, $\Phi_{\left| S\right|,S}$ is the sum over all spanning trees of $\mathcal{G}$ (cf.~Eq.~\eqref{def}). 
Thus, Kirchhoff's Theorem \cite{Kirchhoff}, which has previously been used for the analysis of dynamical systems \cite{Schnakenberg, Li},
can be regarded as a special case of the more general Eq.~\eqref{det_n_2}. 

\begin{figure}
\begin{center}
\includegraphics[width=0.48\textwidth]{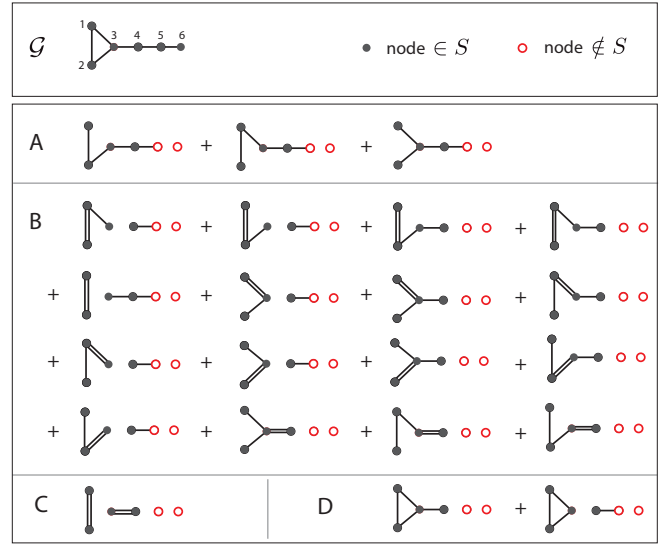}
\end{center}
\vspace{-0.3cm}
\caption{Symbolic calculation of a minor using the zero-row-sum condition. Shown is the graph $\mathcal{G}$, defined by the off-diagonal entries of Eq.~\eqref{ZeroRowSumJacobian}. The terms of the minor $D_{4,S}$ can be written as $\times\cdot\times\cdot\times\cdot\times=A+B+C+2D$, $-\times\cdot\times \cdot\ |=-(B+2C)$, $|\cdot|= C$, $2\times \cdot\ \triangle= -2D$ and $-2\Square= 0$ (cf.~Eq.~\eqref{EqFig2}). It thus follows that $D_{4,S}\equiv\Phi_{4,S}=A$ is the sum over all acyclic subgraphs of $\mathcal{G}$ meeting conditions (i)--(iv).\label{Subgraphs}}
\end{figure}

\vspace{\baselineskip}
\begin{figure}
\begin{center}
\includegraphics[width=0.48\textwidth]{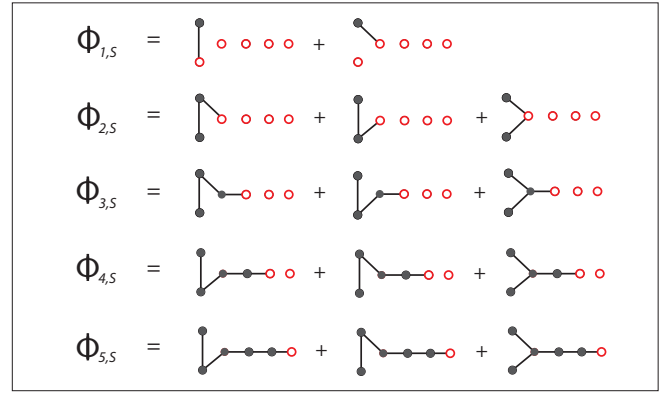}
\end{center}
\vspace{-0.3cm}
\caption{Topological equivalents of minors. Shown is the complete sequence of \mbox{$\Phi_{\left| S\right|,S}$, $\left| S\right|=1,\ldots,5$} for the graph $\mathcal{G}$ from Fig.~\ref{Subgraphs}. \label{Subgraphs2}}
\end{figure}

The simplification of the minors due to the zero-row-sum condition as well as the relation between the $D_{\left| S\right|,S}$ and their topological equivalents $\Phi_{\left| S\right|,S}$ can be illustrated by means of a simple example.
Consider the symmetric zero-row-sum $6\times 6$ matrix
\begin{equation}\label{ZeroRowSumJacobian} 
{\rm \bf J}\!=\!\left(\begin{matrix}
 (J_{11})&J_{12}&J_{13}&0&0&0 \\
 J_{12}&(J_{22})&J_{23}&0&0&0  \\
 J_{13}&J_{23}&(J_{33})&J_{34}&0&0\\ 
0&0&J_{34}&(J_{44})&J_{45}&0\\
 0&0&0&J_{45}&(J_{55})&J_{56}\\
 0&0&0&0&J_{56}&(J_{66})\\
\end{matrix}\right)\ . 
\end{equation}
Here, $(J_{ii})$ abbreviates $-\sum_{j\neq i}J_{ij}$. 
In Fig.~\ref{Subgraphs}, we calculate the minor $D_{4,S=\left\{1,\ldots,4\right\}}$ in terms of the subgraphs of the corresponding graph $\mathcal{G}$.
The calculation illustrates the reasoning that lead to the Eqs.~\eqref{def} and~\eqref{det_n_2}.

On the one hand, the complete sequence of minors $D_{\left| S\right|,S}$, $\left| S\right|=1,\ldots,{\rm rank}({\rm \bf J})$ can be calculated as 
\begin{align}
	D_{1,S_1}&=(-1)\left(J_{12}+J_{13}\right)\nonumber\\
	D_{2,S_2}&=(-1)^2\left(J_{12}J_{13}+J_{12}J_{23}+J_{13}J_{23}\right)\nonumber\\
	D_{3,S_3}&=(-1)^3\left(J_{12}J_{13}+J_{12}J_{23}+J_{13}J_{23}\right)J_{34}\nonumber\\
  D_{4,S_4}&=(-1)^4\left(J_{12}J_{13}+J_{12}J_{23}+J_{13}J_{23}\right)J_{34}J_{45}\nonumber\\
  D_{5,S_5}&=(-1)^5\left(J_{12}J_{13}+J_{12}J_{23}+J_{13}J_{23}\right)J_{34}J_{45}J_{56} \nonumber
\end{align}
where $S_k=\left\{1,\ldots,k\right\}$, and ${\rm rank}({\rm \bf J})=5$ due to the zero-row-sum condition. 

On the other hand, we can use the definition~\eqref{def} to construct the sequence $\Phi_{\left| S\right|,S}$, $\left| S\right|=1,\ldots,{\rm rank}({\rm \bf J})$, directly from the graph $\mathcal{G}$ (cf.~Fig.~\ref{Subgraphs2}). A comparison of both, the algebraic and the topological results, reproduces Eq.~\eqref{det_n_2}.
We note that labeling the nodes in different order would have yielded different algebraic as well as topological expressions.

\section{Topological stability conditions}
Let us shortly summarize what we obtained so far.
The topological reading of determinants maps a Hermitian Jacobian ${\rm \bf J}$ with zero row sum onto a graph $\mathcal{G}$, whose weight matrix is given by the off-diagonal part of $\rm \bf J$.
The minors of ${\rm \bf J}$ can then be interpreted as sums over terms associated with subgraphs of $\mathcal{G}$. 
Combining the Eqs.~\eqref{stability_cond} and~\eqref{det_n_2}, the algebraic stability constraints on the minors of ${\rm \bf J}$ translate into
\begin{equation}\label{stability_cond2}
	\Phi_{\left| S\right|,S}>0, \quad \forall S,  \ \left| S\right|=1,\ldots, {\rm rank}({\rm \bf J}).
\end{equation}
 
We emphasize that the graph $\mathcal{G}$ is not an abstract construction, but combines information about the physical interaction topology and the dynamical units.
For example, if a graph $\mathcal{G}$ has disconnected components, there is a reordering of the variables $x_i$, such that ${\rm \bf J}$ is block diagonal. This implies that the spectra of different topological components of $\mathcal{G}$ decouple and can thus be treated independently.

From Eq.~\eqref{stability_cond2} we can immediately read off a weak \emph{sufficient} condition for stability:
Because $\Phi_{\left| S\right|,S}$ is a sum over products of the $J_{ij}$, a Jacobian with $J_{ij}\geq 0$ $\forall i,j$ is a solution to Eq.~\eqref{stability_cond2} irrespective of the specific structure of $\mathcal{G}$  \cite{remark2}.
By contrast, if $J_{ij}<0$ for some $i,j$, then the existence of solutions of Eq.~\eqref{stability_cond2} is dependent on the topology.

The $\Phi$-notation allows to investigate which combinations of negative $J_{ij}$ in a graph $\mathcal{G}$ lead to the violation of at least one of the Eqs.~\eqref{stability_cond2}: The definition of $\Phi$ implies that any given $\Phi_{\left| S\right|,S}$ can be represented as a function over a set of $\Phi_{\left| S_i\right|,S_i}$
\begin{equation}\label{polynomial}
	\Phi_{\left| S\right|,S}=f\left(\left\{\Phi_{\left| S_i\right|,S_i}\right\}\right),
\end{equation}
where the $S_i$ are subsets of $S$, and $f$ depends on the choice of the $S_i$. 
The key point is now that for showing that a given topological arrangement of negative $J_{ij}$ is incompatible with the stability conditions Eqs.~\eqref{stability_cond2}, it suffices to show that there exists some $\Phi_{\left| S\right|,S}$ and some expansion \eqref{polynomial} thereof, for which there is no solution with all occurring $\Phi$ $>0$.
Let us emphasize that the sketched approach dispenses with the need to specifically calculated any $\Phi_{\left| S\right|,S}$, and hence with addressing the up to $\left| S\right|!$ terms they subsume.  

In \cite{Do}, we employ the outlined approach to show that stability requires that every component of $\mathcal{G}$ has a spanning tree, whose links all have positive weights.
A detailed mathematical proof of this result will be published in \cite{Epperlein}. 
For the sake of completeness we present a sketch of a proof in Appendix A. The proof uses the $\Phi$-notation to show that at least one of the necessary conditions must be violated when no spanning tree of positive connections exists. It thereby also provides an illustration of the direct application of the notation in a calculation.
 
\begin{figure}
\begin{center}
\includegraphics[width=0.47\textwidth]{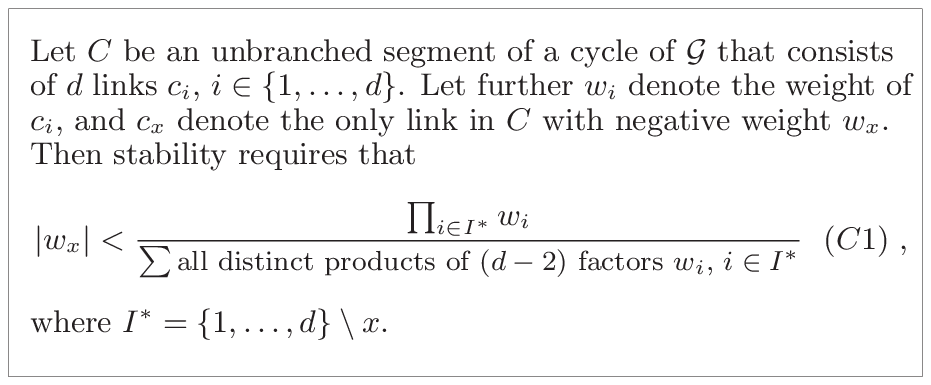}
\end{center}
\sffamily{\caption{Stability sets an upper bound on the absolute value of a negative link weight. The bounding relation is derived in Appendix B.
\label{Upperbound}}}
\end{figure}

In addition to the spanning tree criterion for stability, which pertains to a global property of $\mathcal{G}$, Eq.~\eqref{stability_cond2} implies further restrictions, which pertain to mesoscale properties of $\mathcal{G}$.
Thus, the weights of all links that are not part of any cycle of $\mathcal{G}$ have to be positive.
Moreover, at most one of the links that build an unbranched segment of a cycle of $\mathcal{G}$ may have a negative weight.
Finally, the weight of a `negative link' in an unbranched segment of a cycle of $\mathcal{G}$ is bounded below by a value that depends on the weights of the other links in the segment (cf.~Fig~\ref{Upperbound}). 

Note that while the mesoscale criteria restricting the number and position of negative links can be subsumed under the global spanning tree criterion, the mesoscale criteria restricting the weights of possible negative links are inherently bound to the mesoscale.
In particular, the latter can only be derived by considering $\Phi_{\left| S\right|,S}$ with small and intermediate $\left| S\right|$.
This highlights the benefits of the proposed approach: The JSC provide stability criteria on all scales $\left| S\right|$, which are made accessible by means of the graphical notation and symbolic calculus.

Applied to the Kuramoto model defined in Eq.~\eqref{ODE}, the spanning tree criterion for stability allows to deduce properties shared by all possible phase-locked states \cite{Do}.
In a phase-locked state, $J_{ij}=A_{ij}\cos(x_j-x_i)$, where $x_j-x_i$ is the stationary phase difference between oscillator $j$ and $i$.
Given that all link weights $A_{ij}\geq 0$, stability of the phase-locked state thus requires that the coupling network has a spanning tree of oscillators obeying $|x_j-x_i|<\pi/2$.

\section{General systems of coupled units}

The derivation of stability criteria by graphical interpretation of the JSC is applicable to all systems with Hermitian Jacobian matrix. 
Further, the simplification leading up to the spanning tree criterion is possible, whenever the Jacobian has zero row sums. This condition is satisfied for instance by all systems of the form
	\begin{equation} \label{generalization} \dot{x_i}=C_i+\sum_{j\neq i} A_{ij}\cdot O_{ij}(x_j-x_i)\ , \quad \forall i\in 1\ldots N	\ , 
\end{equation}
where the $A_{ij}$ are the weights of a symmetric interaction network and the $O_{ij}$ odd functions.
We emphasize that the approach remains applicable in heterogeneous networks containing different link strengths $A_{ij}$, coupling functions $O_{ij}$  or intrinsic parameters $C_i$. 

The class of systems to which the present results are directly applicable include general networks of phase oscillators as well as other models such as continuous-time variants of the Deffuant model of opinion formation \cite{Deffuant} or a class of ecological metapopulation models \cite{Leibold}.
 
Although the specific criteria derived above are contingent on the zero-row-sum condition, it can be expected that the general approach proposed here is also applicable to situations where this condition is violated, such as in the model of cooperation among interacting agents studied in \cite{DoRudolfGross}. 
A simple extension of Eq.~\eqref{generalization} which violates the zero-row-sum condition is found by replacing $C_i$ by a function of $x_i$, which allows dynamical retuning of the intrinsic frequency, e.g. for modeling homeostatic feedback in neural networks. 
For illustrating the application of the proposed method to models, in which the zero-row-sum condition is violated in some rows, we consider another model inspired by neuroscience in the subsequent section.

\section{Adaptive Kuramoto model}
We now apply the proposed approach to an example of an adaptive network, i.e., a system in which the topology of the network coevolves with the dynamics of oscillators \cite{ItoKaneko,ZhouKurths,AN}. In the context of the Kuramoto model, adaptive coupling has recently attracted keen interest as it allowed to show that the emergence of synchronous motion can be intimately related to a selection mechanism of specific network topologies \cite{Sendina-Nadal}, and to the identification of complex hierarchical structures in the graph connectivity \cite{ZhouKurths,Almendral}.

We consider a system of $N$ phase oscillators that evolve according to Eq.~\eqref{ODE}, while the coupling strength $A_{ij}$ evolves according to
\begin{equation}\label{Aij_te}
	{\rm \frac{d}{dt}}	A_{ij}=\cos(x_j-x_i)-b\cdot A_{ij}.
\end{equation}
The first term in Eq.~\eqref{Aij_te} states that the more similar the phases of two nodes the stronger reinforced is their connection, the second term guarantees convergence. In a stationary, phase-locked state state $A_{ij}=\cos(x_j-x_i)/b$ and all oscillators oscillate with a common frequency $\Omega=\tfrac{1}{N}\sum_i\omega_i$. The stability of this state is governed by a symmetric Jacobian
\begin{equation}\label{Jacobian}
{\rm \bf J}=\left(\begin{array}{cccccc}
 -b&0&0&\multicolumn{1}{|c}{s_{21}}&s_{12}&0 \\
 0&-b&0&\multicolumn{1}{|c}{s_{31}}&0&s_{13}  \\
 0&0&-b&\multicolumn{1}{|c}{0}&s_{32}&s_{23}\\ \cline{1-6}
 s_{21}&s_{31}&0&\multicolumn{1}{|c}{m_{1}}&o_{12}&o_{13}\\
 s_{12}&0&s_{32}&\multicolumn{1}{|c}{o_{12}}&m_{2}&o_{23}\\
 0&s_{13}&s_{23}&\multicolumn{1}{|c}{o_{13}}&o_{23}&m_{3}\\
\end{array}\right)\ ,
\end{equation}
where $o_{ij}:=\tfrac{1}{b}\cos^2(x_j-x_i)$, $m_{i}:=-\sum_{j\neq i} o_{ij}$, $s_{ji}:=\sin(x_j-x_i)$ and we have chosen $N=3$ for illustration. The marked partitioning separates two blocks on the diagonal. The upper one is a diagonal submatrix of size $L\times L$, $L:=N(N-1)/2$, the lower one is a $N\times N$ symmetric submatrix with zero row sum, which we denote as ${\rm \bf j}$.

Let us start our analysis by focusing on the upper left block of ${\rm \bf J}$. In the chosen ordering of variables,
the first $L$ minors $D_{\left| S\right|,S}$ satisfy the stability condition Eq.~\eqref{stability_cond} iff $b>0$.
Concerning the minors of order $\left| S\right|>L$, the following conventions prove advantageous: 
Below, we consider sets $S$ that contain all variables $A_{ij}$ and $n$ of the variables $x_{i}$. 
For every such set, we define $S'$ as the subset of $S$ that only contains the $n$ variables $x_{i}$. 
Moreover, we define $\tilde{{\rm \bf j}}$ as the matrix that is obtained from ${\rm \bf j}$ if all $o_{ij}$ are substituted by $\cos(2(x_j-x_i))$. 
Finally, we define $\tilde{D}_{n,S'}$ as the determinant of the submatrix of $\tilde{{\rm \bf j}}$, which is spanned by the variables $\in S'$. 

We find that
\begin{equation}
 D_{L+n,S}=(-1)^{L}b^{L-n}\;\cdot \tilde{D}_{n,S'}.
\end{equation}
As $\tilde{{\rm \bf j}}$ is symmetric and has a zero row sum, its minors, $\tilde{D}_{n,S}$, can be rewritten using Eq.~\eqref{det_n_2}
\begin{equation}\label{det_n_3}
	D_{L+n,S}=(-1)^{L+n}b^{L-n} \tilde{\Phi}_{n,S'},
\end{equation}
where $\tilde{\Phi}_{n,S'}$ refers to subgraphs of the graph $\tilde{\mathcal{G}}$ defined by the off-diagonal entries of $\tilde{{\rm \bf j}}$.
Stability requires that ${\rm sgn}\left(D_{L+n,S}\right)={\rm sgn}\left((-1)^{L+n}\right)$. As the necessary stability condition $b>0$ implies  $b^{L-n}>0$, it follows that in a stable system
\begin{equation}\label{stability_cond_adap}
\tilde{\Phi}_{n,S'}>0, \quad \forall S,\  n=1,\ldots ,{\rm rank}(\tilde{{\rm \bf j}}).
\end{equation}
Comparison with Eq.~\eqref{stability_cond2} reveals that a necessary condition for stability is that every component of the graph $\tilde{\mathcal{G}}$ has a positive spanning tree. 
Revisiting the definitions of graph $\tilde{\mathcal{G}}$, we find that the weight of a link $ij$ of $\tilde{\mathcal{G}}$ is given by $\cos(2(x_j-x_i))$.
Hence, every component of $\tilde{\mathcal{G}}$ has a positive spanning tree iff every component of the adaptive coupling network has a spanning tree of oscillators obeying $|x_j-x_i|<\pi/4$.
The restriction on the stationary phase-differences in a stable, phase-locked state are thus more strict in the adaptive than in the non-adaptive case. 

\section{Applicability of the graphical notation to isospectrality problems}

In addition to stability analysis, the graphical notation introduced here allows for instance exploring the isospectrality of Hermitian or symmetric matrices \cite{Eskinetal,Shirokov}, which differ with respect to the signs of some off-diagonal entries.
To see this consider a Hermitian matrix ${\rm\bf A}\in\mathbb{C}^{n\times n}$. Its characteristic polynomial $\chi$ can be calculated as 
	\[\chi(\lambda)=D_n({\rm\bf A}-\lambda {\rm\bf I})
\]
Considering the structure of the graph $\mathcal{G}$ associated to ${\rm\bf A}-\lambda {\rm\bf I}$ allows to determine which symbols contribute to $\chi$. One can then determine which changes of sign of off-diagonal entries leave all contributing symbols and thus the characteristic polynomial and the spectrum invariant.

For instance, if the graph $\mathcal{G}$ is a tree, then the only symbols that contribute to the characteristic polynomial $\chi$ are the symbols $\times$ and $|$. Thus, $\chi$ is a polynomial of factors $(A_{ii}-\lambda)$ (corresponding to symbols $\times$), and factors $A_{ij}A_{ji}$ (corresponding to symbols $|$). Due to the hermiticity of ${\rm\bf A}$, $A_{ij}A_{ji}=|A_{ij}|^2$.
It follows that $\chi$, and therewith the spectrum, is invariant under any operation that changes the sign of a pair of off-diagonal entries $A_{ij}\rightarrow -A_{ij}, \ A_{ji}\rightarrow -A_{ji}$.

Along the same line, we can infer isospectrality relations for matrices ${\rm\bf A}$, whose corresponding graphs $\mathcal{G}$ are composed only of tree-like subgraphs and isolated cycles. For such matrices, the spectrum is invariant under symmetry preserving sign changes of
\begin{itemize}
	\item  an arbitrary number of off-diagonal entries that do not belong to cyclic subgraphs of $\mathcal{G}$.
	\item  an even number of off-diagonal entries that belong to the same cyclic subgraph of $\mathcal{G}$.
\end{itemize}

\section{Conclusions}
In the present paper we presented a general approach for the analysis of symmetrically coupled systems and a specific simplification for
systems, in which the Jacobian matrix ${\rm \bf J}$ additionally has a zero row sum. 

Using a graphical interpretation of Jacobi's signature criterion we showed that for all non-trivial eigenvalues of ${\rm \bf J}$ to be negative, the graph $\mathcal{G}$, whose adjacency matrix is given by the off-diagonal part of ${\rm \bf J}$, has to obey necessary topological conditions. 
Thus, $\mathcal{G}$ must have a spanning tree of links with positive weights, which restricts the number and position of potential negative entries in the Jacobian matrix. 
Moreover, the absolute value of potential negative link weights is bounded above by a topology depended relation. 

We used the approach for deriving necessary conditions for local asymptotic stability of stationary and phase-locked states in networks of phase oscillators. 
Our results provide an analytical angle that is complementary to statistical analysis of network synchronizability. Where statistical approaches reveal global features impinging on the propensity to synchronize, our approach can pinpoint specific defects precluding synchronization. We note that such defects can occur on all scales, corresponding to the violation of the signature criterion in subgraphs of different size. This highlights synchronization of phase oscillators as a simple but intriguing example in which instabilities can arise from local, global or mesoscale structures. In the future, the approach proposed here may provide a basis for further investigation of these instabilities.  

A limitation of the present approach is that the Jacobian matrix of the state under consideration has to be available in some form. 
However, this limitation is inherent to all forms of linear stability and bifurcation analysis and several strategies have been developed for mitigating it. These include for instance the equation-free approach \cite{Kevrekidis}, where the desired information on the Jacobian is generated on-demand from simulations or experiments, and generalized modeling \cite{RudolfGross}, which reveals a parameterized presentation of the Jacobian for very general classes of systems. Even if these specific techniques are not applied, researchers are often aware of the structure of the Jacobian in the system under consideration, which can be sufficient for gaining some fundamental insights by the proposed approach. 

It is apparent that applying the stability analysis proposed here may require overcoming specific problems. 
In particular, the requirements of symmetry and zero row sum may necessitate extending the basic scheme or focusing on specific cases, as we did in the adaptive network example. 
However, if these problems can be overcome, the approach can reveal hard analytical conditions that all stable steady states must meet. Further, it can identify topological structures (on all scales) that prevent the existence of such states and pinpoint them in networks. 
Because this type of information is complementary to insights revealed by common statistical approaches, it potentially of high value for the analysis of the system.

\begin{appendix}
\section{Spanning Tree Criterion for stability}
\begin{figure*}
\begin{center}
\includegraphics[width=0.9\textwidth]{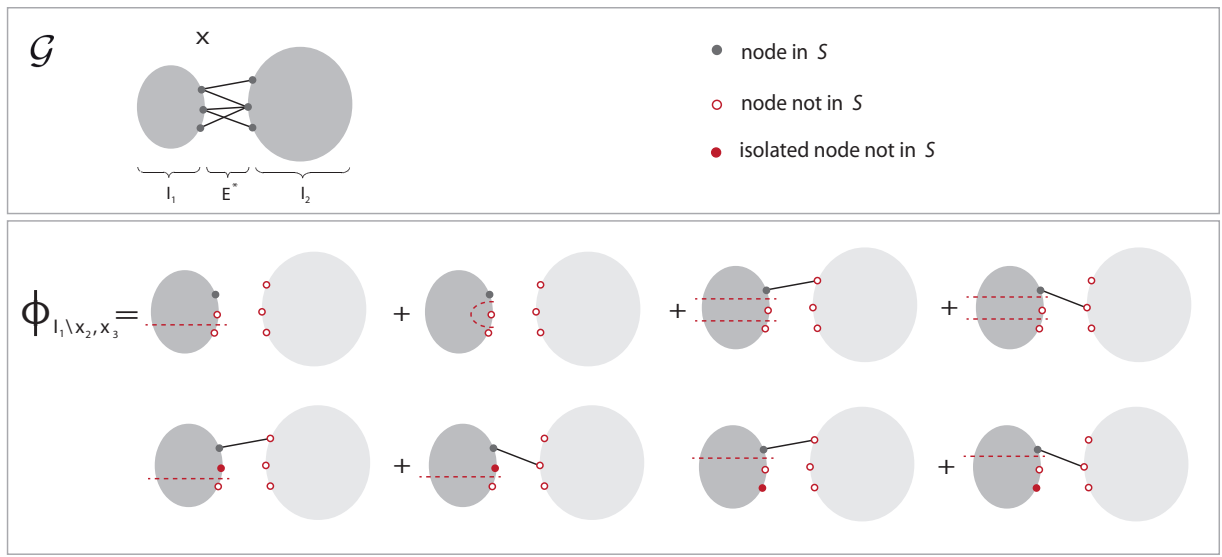}
\end{center}
\caption{Expanding $\Phi_{S}$. Upper panel: Consider a graph $\mathcal{G}$ without positive spanning tree. The grey shaded areas represent the subgraphs of $\mathcal{G}$ induced by $I_1$, $I_2$ respectively. The nodes in $X$ shall be labeled $x_1,x_2,x_3$ with $x_1$ being the upmost and $x_3$ being the bottommost node. 
Lower panel: We now want to symbolically construct all terms that contribute to $\Phi_{I_1\setminus \left\{x_2,x_3\right\}}$. To this end, we use a definition of $\Phi_{S}$, that is equivalent to the one given in the main text:
$\Phi_{S}$ is the sum over all forests of $\mathcal{G}$, in which each tree contains exactly one vertex $\notin S$. 
The panel displays the different realizations of such forests with two trees (term 1--2), three trees respectively (term 3--8).
Here, a grey shaded area divided by $n$ dashed lines symbolizes all possibilities to span $I_1\setminus \left\{x_2,x_3\right\}$ with $n$ trees such that no two nodes in $X$, which are divided by a dashed line, are part of the same tree.  
Comparing with the definition of $\tau_Y$ and $\sigma_Y$ reveals that the sum over the first two symbolic terms equals the $B=\emptyset$ term of Eq.~\eqref{expansion}, while the sum over all other symbolic terms equals the $B=\left\{x_1\right\}$ term of Eq.~\eqref{expansion}.    
\label{Expansion_Fig}}
\end{figure*}

To prove the spanning tree criterion for stability consider that in a network without a positive spanning tree it must be possible to partition the nodes into two nonempty sets $I_1$, $I_2$ such that  
\begin{equation} \label{negativecut}
	J_{ij}\leq 0 \quad \forall \ i,j \ | \  i\in I_1,\  j\in I_2  \ .
\end{equation}
The idea is now to evaluate the stability conditions~\eqref{stability_cond2} for different $S\subseteq I_1$ thereby exploiting that all links leading out of $I_1$ have negative weights. 
For this purpose,it is convenient to define $E^*$ as the set of links connecting $I_1$ and $I_2$, and $X=\left\{x_{1},\ldots,x_m\right\}$ as the subset of nodes $\in I_1$ incident to at least one link from $E^*$ (`boundary vertices'). 
Further, we define $\sigma_i$ as the sum over all elements of $E^*$ incident to $x_i$, and, for any subset $Y$ of $X$, $\sigma_{Y}:=\prod_{m\in Y}\sigma_m$. 
Finally, for any subset $Y$ of $X$, we define $\tau_{Y}$ as the sum over all forests of $\mathcal{G}$ that (i) span $I_1$, and (ii) consist of $\left|Y\right|$ trees each of which contains exactly one element from $Y$. 

With the above definitions we can write
\begin{equation}\label{expansion}
	\Phi_{I_1\setminus C}=\sum_{B\subseteq X\setminus C}\sigma_B\tau_{B\cup C} \ ,
\end{equation} 
where $B$ and $C$ are disjoint subsets of $X$ and $\Phi_{I_1\setminus C}:=\Phi_{q,I_1\setminus C}$. The rather abstract equation is illustrated in Fig.~\ref{Expansion_Fig} by means of an example. 

Equation~\eqref{expansion} can be read as an expansion of $\Phi_{I_1\setminus C}$ in contributions from $E^*$. 
This expansion has the advantage that the sign of one of the two factors in each term is known: ${\rm sgn}\left(\sigma_B\right)=(-1)^{|B|}$. 
The signs of the factors $\tau_{B\cup C}$, however, are still indeterminate.
Nevertheless, we can show that Eq.~\eqref{expansion} is incompatible with the stability condition \eqref{stability_cond2} by considering a linear combination of $\Phi_{I_1\setminus C}$:
\begin{align}
	\sum_{C\subseteq X}\underbrace{(-1)^{\left|C\right|}\sigma_{C}}_{\substack{>0\ \text{by} \\\text{construction}}}\underbrace{\Phi_{I_1\setminus C}}_{>0}&=\sum_{C\subseteq X}(-1)^{\left|C\right|}\sigma_{C}\sum_{B\subseteq X\setminus C}\sigma_B\tau_{B\cup C} \nonumber\\
	&=\sum_{\substack{C\subseteq X\\B\subseteq X\setminus C}} (-1)^{\left|C\right|}\sigma_{B\cup C} \tau_{B\cup C}\nonumber\\
	&=\sum_{\substack{A\subseteq X\\C\subseteq A}} (-1)^{\left|C\right|}\sigma_{A} \tau_{A}\nonumber\\
	&=\sum_{A\subseteq X}\sigma_{A} \tau_{A}\underbrace{\sum_{C\subseteq A}(-1)^{\left|C\right|}}_{=0}=0 \ , \nonumber
\end{align}
which is a contradiction.
Therefore the existence of a spanning tree of positive elements is a necessary condition for stability. 

\section{Lower bound on the negative link weights}
Consider a path of $d-1$ degree-two nodes $v_i$, $i\in \left\{1,\ldots,d-1\right\}$, that are part of at least one cycle of $\mathcal{G}$. 
Together with the $d$ egdes $c_i$ that are incident to at least one of the nodes $v_i$, the path constitutes an unbranched segment $C$ of a cycle of $\mathcal{G}$ (cf.~Fig.~3).
 
According to the spanning tree criterion, $C$ can have at most one link with negative weight. 
Below, we consider the case where $C$ has exactly one such link, and show that Eq.~\eqref{stability_cond2} imposes an upper bound on the absolute value of the negative link weight.

We use the following conventions: 
Let the nodes be labeled such that the indices $i$ occur in an increasing order if $C$ is paced out. 
And let the links be labeled such that $v_i$ is incident to $c_i$ and $c_{i+1}$ (cf.~Fig.~2S). 
Further, let $w_i$ denote the weight of $c_{i}$. And lastly, let $c_x$ denote the only link in $C$ with negative weight $w_x$. 
Below, we show that stability requires that
	\begin{equation}\label{upperbound}\textstyle{	
	\left|w_{x}\right|<\frac{\prod \text{all $w_i$, $i\in I^*$} }{\sum \text{all distinct products of $(d-2)$ factors $w_i$, $i\in I^*$}} \ , }
\end{equation}
where $I^*= \left\{1,\ldots,d\right\}\setminus x$.

For deriving Eq.~\eqref{upperbound}, we consider a sequence of conditions $\Phi_{S_1}\ldots\Phi_{S_{d-1}}$. The sequence is constructed as follows: we choose $S_1=x$, $S_2=\left\{x,x+1\right\}$, $S_3=\left\{x,x+1,x+2\right\}$, and so forth until $S_{d-x}=\left\{x,x+1,x+2,\ldots,d-1\right\}$. The remaining elements of the sequence are then constructed as $ S_{d-x+1}=S_{d-x}\cup \left\{x-1\right\}$, $S_{d-x+2}= S_{d-x}\cup\left\{x-1,x-2\right\}$, and so forth until $ S_{d-1}=\left\{v_i\right\}_{i\in \left\{1,\ldots,d-1\right\}}$.

The first element of the sequence, $\Phi_{S_1}>0$, stipulates that $w_x+w_{x+1}>0$, and thus that $-w_x<w_{x+1}$. 

The second element of the sequence $\Phi_{S_2}>0$ stipulates that $w_xw_{x+1}+w_xw_{x+2}+w_{x+1}w_{x+2}>0$ and thus that $-w_{x}<(w_{x+1}w_{x+2})/(w_{x+1}+w_{x+2})$.

More generally, every element $\Phi_{S_i}$ of the sequence $\Phi_{S_1},\ldots,\Phi_{S_{d-1}}$ modifies the upper bound on $-w_x$ as per
\begin{equation}\label{upperbound2}\textstyle{
	-w_{x}<\frac{\prod \text{all $w_i$, $i\in S_i^*$} }{\sum \text{all distinct products of $(\left|S_i\right|-1)$ factors $w_i$, $i\in S_i^*$}} \ , }
\end{equation} 
where $S_{i}^*=(S_i\setminus x)\cup(\max(S_i)+1)$.

The right hand side of Eq.~\eqref{upperbound2} is monotonously decreasing with increasing $\left|S_i\right|$. We can thus conclude that $C$ can have at most one link with negative weight, whose absolute value $\left|w_x\right|$ is bounded above by Eq.~\eqref{upperbound}.

\end{appendix}


\end{document}